\documentclass[a4paper,12pt]{article}
\usepackage{moreverb}
\usepackage{lineno}
\usepackage{graphicx}
\usepackage[english]{babel}
\usepackage[latin1]{inputenc}
\usepackage{amsmath}
\usepackage{amssymb}
\usepackage{multirow}
\usepackage{amsthm}
\usepackage{amsfonts}
\usepackage{graphicx}
\theoremstyle{plain}
\newtheorem{theorem}{Theorem}[section]
\newtheorem{definition}{Definition}[section]
\newtheorem{lemma}[theorem]{Lemma}

\title{Markov Chain Order  Estimation and $\chi^2-divergence$ measure}
\author{ A.R. Baigorri\thanks{baig@unb.br}\\ Mathematics Department \\ UnB 
        \and C.R. Gon\c{c}alves \thanks{The author was partially supported by {\scriptsize PROCAD/CAPES}, {\scriptsize CASADINHO/CAPES} 
         and {\scriptsize PRONEX/FAPDF} --- catia@mat.unb.br}
             \\ Mathematics Department \\ UnB 
        \and P.A.A. Resende \thanks{pa@mat.unb.br} \\ Mathematics Department \\ UnB }
\date{March 01, 2012}

\begin{document}
\maketitle
\linenumbers

\begin{abstract}
We use the $\chi^2-divergence$  as a measure of diversity between probability 
densities  and review the basic properties of the estimator $\Delta_2(.\Vert.).$ 
We define a few objects which capture relevant information from the sample of a Markov Chain to be used 
in the definition of a couple of  estimators i.e.
the \textit{Local Dependency Level} and \textit{Global Dependency Level} for a Makov chain sample. After exploring their
properties  we propose  a new estimator for the  Markov chain order. Finally we show a few tables containing numerical 
simulation results, comparing the perfomance of the new estimator with the well known and already established AIC, BIC and EDC estimators.
\end{abstract}

\vspace{20pt}
\section{Introduction}
A Markov Chain is a discrete stochastic process $\mathbb{X}= {\{ X_n \}}_{n \geq 0} $ with   state space $E$,
cardinality $\vert E \vert < \infty$ for which there is a $ k \geq 1 $ such that for $n\geq k, \,\,(x_1,....,x_n) \in E^n$ 

\begin{linenomath*}
\begin{equation*}
P(\textbf{X}_1 = x_1,..,\textbf{X}_n = x_n)= P( \textbf{X}_1 = x_1,..,\textbf{X}_k = x_k)
{\Pi}_{i=k+1}^n Q(x_i \vert x_{i-k},..., x_{i-1})
\label{def_order}
\end{equation*}
\end{linenomath*}

for suitable transition probabilities $ Q(. \vert .)$. The class of processes that holds the above condition for a given
$k \geq 1$ will be denoted by ${\cal{M}}_k$, and ${\cal{M}}_0$ will denote the class of i.i.d. processes.
The order of a process in $\cup_{i=0}^\infty \,{\cal{M}}_i$ is the smallest integer $\kappa$ such that 
$\mathbb{X}= {\{ X_n \}}_{n \geq 0} \in {\cal{M}}_\kappa$.

\vspace{10pt}
Along the last few decades there has been a great number of research on the estimation of the order 
of a Markov Chains,
starting with  M.S. Bartlett \cite{Bartlett}, P.G. Hoel \cite{Hoel}, I.J. Good \cite{Good}, 
T.W. Anderson \& L.A. Goodman \cite{Anderson-Goodman}, P. Billingsley \cite{Billingsley1},
\cite{Billingsley2} among others, and more recently, H. Tong \cite{Tong}, 
G. Schwarz \cite{Schwarz}, R.W. Katz \cite{Katz}, I. Csiszar and P. Shields \cite{CsiszarShields_1}, 
L.C. Zhao et all \cite{Zhao} had contributed with new Markov chain order estimators.

Since 1973, H. Akaike \cite{Akaike}  entropic information criterion, known as AIC, has had a fundamental impact
in statistical model evaluation problems. The AIC has been applied by Tong, for example, to the 
problem of estimating the order of autoregressive processes, autoregressive integrated moving average 
processes, and Markov chains. The Akaike-Tong (AIC) estimator was  derived as an asymptotic approximate estimate of the 
Kullback-Leibler information discrepancy and provides a useful tool for evaluating models estimated by 
the maximum likelihood method. Later on,  Katz derived the asymptotic distribution of the 
estimator and showed  its inconsistency, proving that there is a positive probability of overestimating 
the true order no matter how large the sample size. 
Nevertheless, AIC is the most used and succesfull Markov chain  order estimator used at the present time, 
mainly because it is more  efficient than BIC  for small sample. 

The main consistent estimator alternative, the BIC estimator, does not perform too well for relatively 
small samples, 
as it was pointed out by Katz \cite{Katz} and Csiszar \& Shields \cite{CsiszarShields_1}. It is natural  to 
admit that the expansion of the Markov Chain complexity (size of the state space and order)  has significant 
influence on the sample size required for the identification of the unknown order, even though, most 
of the time it is difficult  to obtain  sufficiently large  samples.

In this notes we'll use  a different entropic object called $\chi^2-divergence$, 
and  study its behaviour when  applied to samples from random variables with multinomial empirical distributions 
\begin{linenomath*}
$${\cal X} = \{ X_{i} \}_{1 \leq i \leq r}$$ 
\end{linenomath*}
derived from  a Markov Chain sample. 
Finally, we shall propose a new strongly consistent Markov Chain order estimator more efficacious than the
already established AIC and BIC,  which it shall be exhibited through the outcomes of several numerical simulations.

In Section~2 we succinctly review the concept of $f-divergence$ and its properties. In Section~3,  
the $\chi^2$-divergence estimator is defined
reviewing some results concerning its convergence, as well as we briefly elaborate  about the Law of Iterated Logarithm (LIL)
for our particular situation. In Section~4 the Makov chain sample is brought to attention, some notation introduced 
and the estimators \textit{Local Dependency Level}  and \textit{Global Dependency Level}, 
which are the groundsill of the consistent Markov  chain order estimator, subsequently defined. 
Finally, in Section~4 we describe the procedures used and the results obtained in an exploratory
numerical simulations.

\vspace{20pt}
\section{Entropy and f-divergences}

\subsection{Definitions and Notations}
An $f-divergence$ is a function that measures the  discrepancy between two probability distributions $P$ and $Q$. 
The divergence is intuitively an average of the function $f$ of the odds ratio given by $P$ and $Q$.

These divergences were introduced and studied independently by Csiszar, Csiszar\&Shields and 
Ali\&Silvey among others (\cite{Csiszar}, \cite{CsiszarShields}, \cite{AliSilvey}) 
and sometimes are 
referred as \textit{Ali-Silvey distances}.

\vspace{20pt}

\begin{definition}\label{chi_square_empirical_def}
\begin{linenomath*}
Let $P$ and $Q$  be  discrete probability densities  with support $S(P) = S(Q)= E = \{1,...m\}$. 
For $f(t)$  convex function defined for $t~>~0, f(1) = 0$, the $f-divergence$ for the distributions $P$ and $Q$ 
is 
$$ D_f (P \Vert Q) = \sum_{a \in A} Q(a) f\left ( \frac {P(a)}{Q(a)} \right ).$$
Here we take $ 0f(\frac{0}{0}) = 0$ , $ f(0)=\lim_{t \to 0}f(t)$, $  0f(\frac{a}{0}) = \lim_{ t \to 0 } t f(\frac{a}{t})=
a \lim_{ u \to \infty }\frac{f(u)}{u}. {\quad \scriptstyle \blacklozenge}$
\end{linenomath*}
\end{definition}

\vspace{10pt}
For example:
\begin{linenomath*}
$$ f(t) = t \log(t) \Rightarrow D_f(P \Vert Q) = D(P \Vert Q) = \sum_{a \in A} P(a) \log \left ( \frac{P(a)}{Q(a)} \right ),$$
\end{linenomath*}

\begin{linenomath*}
$$ f(t) = (1 - t^2)  \Rightarrow D_f(P \Vert Q) = \sum_{a \in A}  \frac{(P(a) - Q(a))^2}{Q(a)},$$
\end{linenomath*}

which are called \textit{relative  entropy}  and $\chi^2-divergence$, respectively.
From now on the $\chi^2-divergence$ shall be denote by  $ D_2(P \Vert Q)$.

Observe that the triangular inequality is not satisfied in general, so that  $D_2(P \Vert Q)$ defines no distance in the strict sense. 

\vspace{20pt}
A basic theorem  about \textit{f-divergences} is the following approximation by the $D_2(P \Vert Q)$.

\begin{theorem}(Csiszar \& Shields \cite{CsiszarShields})
If f is twice differentiable at t=1 and $f^{''}(1) > 0$ then for any $Q$ with support $S(Q) = A$ and $P$ \textbf{close} to $Q$ 
\begin{linenomath*}
$$ D_f (P \Vert Q) \sim \frac{f^{''} (1)}{2} D_2(P \Vert Q).$$
\end{linenomath*}
Formally, $D_f(P \Vert Q)/D_2(P \Vert Q) \rightarrow f^{''}(1)/2$ as $P\stackrel{D}\rightarrow Q$  ${\quad \scriptstyle \blacklozenge}$
\end{theorem}

The $\chi^2$-square  divergence $D_2(P \Vert Q)$  test is well known  statistical test procedure
close related to the chi-square distribution. See \cite{Pardo} for thorough and detailed references.

\vspace{30pt}
\section{Derived Markov Chains}
\vspace{10pt}
Let $\textbf{X}_1^n = (X_1,...,X_n)$  be a sample from a multiple stationary Markov chain $\mathbb{X}= {\{ X_n \}_{n \geq 1}}$ 
of unknown order ${\kappa}.$  
Assume that  $\mathbb{X}$ take values on a finite state space $E  = \{ 1,2,...,m\}$ with  transition probabilities  given by

\begin{linenomath*}
\begin{equation}
 p({x_{\kappa+1}\vert x_1^\kappa}) = P(X_{n+1}=x_{n+1} \vert X_{n-\kappa+1}^n=x_1^\kappa) > 0
\label{chain_cond_1}
\end{equation}
\end{linenomath*}

where $x_1^\kappa = x_1^j\,x_{j+1}^\kappa = (x_1,...,x_\kappa) \in E^\kappa$.

\vspace{25pt}
Following Doob \cite{Doob}, from the process $\mathbb{X}$ we can derive a first order MC, 
$\mathbb{Y}^{(\kappa)}=\{{Y}^{(\kappa)}_n\}_{n \geq 0}$ by setting 
$ {Y}^{(\kappa)}_n=({X}_n, ....,{X}_{n+\kappa-1})$ so that for $v=(i_1,.....i_\kappa)$ and 
$w=({i^{'}}_1,......,{i^{'}}_\kappa)$ 

\begin{linenomath*}
\begin{equation*}
P(Y_{n+1}^{(\kappa)} = w \vert Y_n^{(\kappa)} = v ) = {\tilde{p}}_{vw} = \begin{cases} p({{i^{'}_\kappa}} \vert i_1....i_\kappa) , 
\,\, {i^{'}_j}=i_{j+1},\,\,\,  j=1,...,(\kappa -1)    \\ 0, 
\,\,\,\, otherwise. \end{cases} 
\end{equation*}
\end{linenomath*}
  
\vspace{10pt}
Clearly $\mathbb{Y}^{(\kappa)}$ is a first order and homogeneous MC that from now on shall be called the derived process, which 
by (\ref{chain_cond_1}) is irreducible and positive recurrent MC having unique stationary distribution, 
say $ \varPi_\kappa $.  It is well known, see [\cite{Doob}-Chap. 5.3], that the derived Markov Chains $\mathbb{Y}^{(l)}, \,\, l \geq \kappa$
is irreducible and aperiodic, consequently ergodic.

There exists an equilibrium (stationary) distribution $\varPi_\kappa(.)$ satisfying for any initial distribution $\nu$ on $E^\kappa$

\begin{linenomath*}
\begin{equation*}
\lim_{n \to \infty} \vert P_\nu( Y_n^{(\kappa)}= x_1^\kappa ) - \varPi_\kappa(x_1^\kappa)\vert = 0, 
\end{equation*}
\end{linenomath*}

and

\begin{linenomath*}
\begin{equation*}
\varPi_\kappa(x_1^\kappa) = \sum_{z_1^\kappa} \varPi_\kappa (z_1^{\kappa}) \, p(x_\kappa \vert z_1^{\kappa})=
\sum_{x} \varPi_\kappa (x\,x_1^{\kappa-1}) \, p(x_\kappa \vert x\,x_1^{\kappa-1}).
\end{equation*}
\end{linenomath*}
 
Likewise, for $\mathbb{Y}^{(l)}, \,\, l > \kappa$

\begin{linenomath*}
\begin{equation}
\varPi_l(x_1^l) =  \varPi_\kappa (x_1^\kappa) \, p(x_{\kappa+1} \vert x_1^{\kappa})...p(x_{l} \vert x_{l-\kappa}^{l-1})
                =  \sum_{x} \varPi_l (x\,x_1^{l-1}) \, p(x_l \vert x\,x_{l-\kappa}^{l-1}).
\label{ergodic_measure_l}
\end{equation}
\end{linenomath*}

which shows that $\varPi_l$ defined above, is a stationary distribution for $\mathbb{Y}^{(l)}$. For the sake of notation's
simplicity we'll use, from now on 
 
\begin{linenomath*}
\begin{equation}
\boxed{\qquad
\varPi(a_1^l) = \varPi_l(a_1^l), \,\,\,\, l \geq \kappa.\qquad}
\label{stationary_dist_notation}
\end{equation}
\end{linenomath*}

\vspace{20pt}
Now, let us return to $\textbf{X}_1^n= (X_1, X_2,...,X_n)$ and define

\begin{linenomath*}
\begin{equation}
N(x_1^l \vert \textbf{X}_1^n) = \sum_{j=1}^{n-l+1} 1(X_j=x_1,...,X_{j+l-1}=x_l)
\label{number_of_strings}
\end{equation}
\end{linenomath*}

i.e. the number of ocurrences of $x_1^l$ in $X_1^n$. If $l=0$ we take $N(\,.\,\vert \textbf{X}_1^n) = n$. The sums are taken over positive terms
$N(x_1^{l+1} \vert \textbf{X}_1^n) > 0$, or else, we convention $0/0$ or $0.\infty$ as $0$.

\vspace{25pt}
Now we define the  empirical random variables $\textbf{X}_{i\,\alpha}$, for $i \in E$ and $\alpha \in E^\eta.$
\vspace{5pt}

\begin{definition}\label{def_X_ialpha}
\begin{linenomath*}
For $\alpha = (a_1,...,a_\eta)= a_1^\eta \in E^\eta$ and $i \in E$, let ${X}_{i\,\alpha}$ be the random variable
taking values in $E$,
extracted from the MC sample $\textbf{X}_1^n$, defined as 
\end{linenomath*}

\begin{linenomath*}
\begin{equation}
P( {X}_{i\,\alpha} = l )= \frac{N(i\,a_1^\eta\,l \, \vert \, \textbf{X}_1^n)}{N(i\,a_1^\eta\,\vert\,\textbf{X}_1^n)}, 
\,\, l \in E.
\label{probability_X_ialpha}
\end{equation}
\end{linenomath*}

with  

$\textbf{X}_{i\,\alpha}=\left ( X_{i \, \alpha}^{(1)},...,X_{i \, \alpha}^{(n_{i \alpha})} \right )$ its sample of 
size $n_{i \, \alpha}$.
\end{definition}

\vspace{10pt}
Observe that for $i,j \in E$

\begin{linenomath*}
\begin{equation*}
\textbf{O}^\alpha_n(i,j)  = N(i\,a_1^\eta\,\vert\,\textbf{X}_1^n)
\end{equation*}
\end{linenomath*}

where $\textbf{O}_n^\alpha$ is the empirical random variables that describe the
$X_{i \, \alpha}, \, 1~\leq~i\leq~m$ observed frequencies. Likewise, we define the expected frequencies

\begin{linenomath*}
\begin{equation*}
\textbf{E}_n^\alpha(i,j) = \frac{ \sum_{l} O_n^\alpha(i,l) \sum_{l} O_n^\alpha(l,j)}{\sum_{kl} O_n^\alpha(k,l)} 
\end{equation*}
\end{linenomath*}

and the respective probability functions

\begin{linenomath*}
\begin{equation*}
\textbf{\textit{P}}_{\textbf{O}_n^\alpha}(i,j) = \frac{\textbf{O}_n^\alpha(i,j)}{N(a_1^\eta\,\vert\,\textbf{X}_1^n)},\,\,\, i,j \in E
\end{equation*}
\end{linenomath*}

\begin{linenomath*}
\begin{equation*}
\textbf{\textit{P}}_{\textbf{E}_n^\alpha}(i,j) = \frac{\textbf{E}_n^\alpha(i,j)}{N(a_1^\eta\,\vert\,\textbf{X}_1^n)},\,\,\, i,j \in E.
\end{equation*}
\end{linenomath*}

\vspace{10pt}
Finally the $\chi^2$-square  divergence 

\begin{linenomath*}
\begin{eqnarray} 
{\hat{\varDelta}}_2(\textbf{\textit{P}}_{\textbf{O}^{\alpha}_n} \Vert \textbf{\textit{P}}_{\textbf{E}^{\alpha}_n}) 
&=&\,n\,\sum_{i=1}^r \sum\limits_{j=1}^m   \frac{ (\textbf{\textit{P}}_{\textbf{O}^{\alpha}_n}(i,j) - 
                                               \textbf{\textit{P}}_{\textbf{E}^{\alpha}_n}(i,j))^2 }
                                              {\textbf{\textit{P}}_{\textbf{E}^{\alpha}_n}(i,j) } \nonumber \\
&=&\,n \,\, \varDelta_2(\textbf{\textit{P}}_{\textbf{O}^{\alpha}_n} \,\Vert \,\textbf{\textit{P}}_{\textbf{E}^{\alpha}_n}).
\label{Delta_2_alpha_estimator} 
\end{eqnarray}
\end{linenomath*}

\vspace{25pt}
Now we  derive a version of the  Law of Iterated Logarithm, significant
for the establisment of subsequent results about the convergence of
${\hat{\varDelta}}_2(\textbf{\textit{P}}_{\textbf{O}^{\alpha}_n} \Vert \textbf{\textit{P}}_{\textbf{E}^{\alpha}_n})$.

\vspace{15pt}
\begin{lemma}\label{ILL_lemma}\cite{MeynTweedie}(Theorems 17.0.1 \& 17.2.2)\,\, Let $ \mathbb{X}=\{X_n\}_{n>0}$ be a ergodic Markov 
chain with finite state space $E$ and stationary distribution 
$\varPi$, $g:~E \longrightarrow \mathbb{R}$, $S_n(g) = \sum_{j=1}^n g(X_j)$ 
and 

\begin{linenomath*}
\begin{equation*}
\sigma_g^2 = E_\pi\,(g^2(X_1)) + 2 \sum_{j=2}^n E_\pi\,(g(X_1) g(X_j)))
\end{equation*}
\end{linenomath*}

then:

\vspace{10pt}
\quad (a) If $\sigma^2_g = 0$, then a.s.  $ \lim_{n \to \infty} \frac{1}{\sqrt{n}} [ S_n(g) - E_\pi (S_n(g)) ] = 0$.

\quad (b) If $\sigma^2_g > 0$, then a.s. 

\begin{linenomath*}
\begin{equation*}
\limsup_{n \to \infty} \frac{ S_n(g) - E_\pi (S_n(g)) }{\sqrt{2 \, \sigma_g^2 \, n \, log(log(n)) }}=1
\end{equation*}
\end{linenomath*}

 and 

\begin{linenomath*}
\begin{equation*}
\liminf_{n \to \infty} \frac{ S_n(g) - E_\pi (S_n(g)) }{\sqrt{2 \, \sigma_g^2 \, n \,log(log(n)) }}=-1,
\end{equation*}
\end{linenomath*}

($E_\varPi$ \,:\, expectation with initial distribution $\varPi$; a.s. \,:\, almost surely).{$\quad \scriptstyle \blacklozenge$}
\end{lemma}

\vspace{15pt}
\begin{lemma}\label{ILL_lemma_1}\cite{Dorea}(Lemma 2)\,\, If $ \mathbb{Y}^{(\kappa)}$ is ergodic then for $\, \eta  \geq \kappa - 1$,
$\alpha~=~ a_1^\eta $  and  \,\, ${ i\,\alpha\,j }=(i,a_1,...,a_\eta,j)= i\,a_1^\eta\,j \in E^{\eta+2}$  we have a.s. 

\begin{linenomath*}
\begin{equation*}
\limsup_{n \to \infty} 
\frac{\big ( N(i\,a_1^\eta\,j\,\vert\,\textbf{X}_1^n) - 
             N(i\,a_i^\eta\,\vert\,\textbf{X}_1^n)\,p(j \big\vert {i\,a_1^\eta} ) \big )^2}{n \log(\log(n)) } =
             2 \, \varPi( {i\,a_1^\eta\,j} )(1-p(j \big\vert {i\,a_1^\eta} )).{\quad \scriptstyle \blacklozenge}
\end{equation*}
\end{linenomath*}
\end{lemma}

\vspace{35pt}
\begin{theorem}\label{Div_conv_null_hypo_theo1}
Let us refer to (\ref{Delta_2_alpha_estimator}) for 
the definition of $\hat{\varDelta}_2 (P_{\textbf{O}_n^\alpha} \Vert P_{\textbf{E}_n^\alpha} )$, 
as well as the beginning of the present section for complementary definitions and references related to the following result:

\vspace{10pt}
If $\kappa \leq \eta$, there exist $ \cal{L} < \infty$ so that for every $\alpha = i_1^\eta \in E^\eta$ 

\begin{linenomath*}
\begin{equation}
\boxed{
P \left ( \limsup_{n \to \infty} \left [ 
\frac{ \hat{\varDelta}_2 (P_{\textbf{O}_n^\alpha} \Vert P_{\textbf{E}_n^\alpha} )}{2\,\log(\log(n))} \right ] \leq \cal{L} \right ) = 1.
\label{hypo_true_conv}
}
\end{equation}
\end{linenomath*}

\vspace{15pt}
If $\eta = \kappa-1 $, there exist $a_1^\eta \,\, {\scriptstyle \&} \,\, i, j , k\neq i$ such that, 
$p(j\,\vert\,i\, a_1^\eta) \neq p(j\,\vert\,k\,a_1^\eta) $, consequently 

\begin{linenomath*}
\begin{equation}
\boxed{
P \left ( \limsup_{n \to \infty} \left [ 
\frac{ \hat{\varDelta}_2 (P_{\textbf{O}_n^\alpha} \Vert P_{\textbf{E}_n^\alpha} )} {2\, log(log(n))} \right ] = \infty \right ) = 1.}
{\quad \scriptstyle \blacklozenge} \label{hypo_false_conv}
\end{equation}
\end{linenomath*}
\end{theorem}

\vspace{15pt}
\textbf{\textit{Proof}:} The following proof shall be divided in the next two cases.

\vspace{15pt}
\textit{\textbf{Case I:}} \,\, $0 \leq \kappa \leq \eta.$

\vspace{5pt}
From (\cite{Zhao}, \textit{Lemma 3.1}) and by  Definition~\ref{divergence_def1} we can calculate

\begin{linenomath*}
\begin{equation*}
\textbf{O}_n^\alpha(i,j)-\textbf{E}_n^\alpha(i,j) = N(i\,a_1^\eta\,j\,\vert\,\textbf{X}_1^n) - 
  \frac{{N(i\,a_1^\eta\,\vert\,\textbf{X}_1^n)}N(a_1^\eta\,j\,\vert\,\textbf{X}_1^n)}{N(a_1^\eta\,\vert\,\textbf{X}_1^n)}
\end{equation*}
\end{linenomath*}

or, in the limit

\begin{linenomath*}
\begin{equation*}
 \lim_{n \to \infty} \Big ( \textbf{O}_n^\alpha(i,j) - \textbf{E}_n^\alpha(i,j) \Big )^2 = 
 \lim_{n \to \infty} \Big ( N(i\,a_1^\eta\,j\,\vert\,\textbf{X}_1^n) - 
                            N(i\,a_1^\eta\,\vert\,\textbf{X}_1^n) \, p(j \big\vert {i\,a_1^\eta}) \Big )^2
\end{equation*}
\end{linenomath*}

\begin{linenomath*}
\begin{eqnarray*}
&&\limsup_{n \to \infty} 
\frac{\Big ( \textbf{O}_n^\alpha(i,j) - \textbf{E}_n^\alpha(i,j) \Big )^2}{n\,\log(log(n))\,\textbf{P}_{\textbf{E}_n^\alpha}(i,j)} = \\
&=&\limsup_{n \to \infty} 
\left [ \frac{\Big ( N(i\,a_1^\eta\,j\,\vert\,\textbf{X}_1^n) - 
N(i\,a_1^\eta\,\vert\,\textbf{X}_1^n)\,p(j \big\vert {i\,a_1^\eta} \Big )^2}{n\,\log(log(n))} 
\frac{1}{\textbf{P}_{\textbf{E}_n^\alpha}(i,j)} \right ].\\
\end{eqnarray*}
\end{linenomath*}

Similarly

\begin{linenomath*}
\begin{eqnarray*}
&&\lim_{n \to \infty}\textbf{P}_{\textbf{E}_n^\alpha}(i,j) = \lim_{n \to \infty} \frac{{\textbf{E}_n^\alpha}(i,j)}{N(a_1^\eta\vert\,\textbf{X}_1^n)}=
\lim_{n \to \infty}\Big( \frac{N(i\,a_1^\eta\vert\,\textbf{X}_1^n)}{N(a_1^\eta\,\vert\,\textbf{X}_1^n)}
\frac{N(a_1^\eta\,j\,\vert\,\textbf{X}_1^n)}{N(a_1^\eta\,\vert\,\textbf{X}_1^n)} \Big)=\\
&=&\lim_{n \to \infty} \Big ( \frac{N(i\,a_1^\eta\,\vert\,\textbf{X}_1^n)}{n}\frac{n}{N(a_1^\eta\,\vert\,\textbf{X}_1^n)}
\frac{N(a_1^\eta\,j\,\vert\,\textbf{X}_1^n)}{N(a_1^\eta\,\vert\,\textbf{X}_1^n)} \Big )=
\varPi(i\,a_1^\eta) \frac{1}{\varPi(a_1^\eta)} \, p(j\,\vert\,a_1^\eta)=\\
&=& \theta(i,j) > 0.
\end{eqnarray*}
\end{linenomath*}

\vspace{10pt}
By (\ref{chain_cond_1}) and Lemma~\ref{ILL_lemma_1} we have that $\,\, \min_{i,j} \theta(i,j) > 0\,\,$ with

\begin{linenomath*}
\begin{equation*}
{\cal L} = \min_{i,j} \theta(i,j)  \sum_{i=1}^m\sum_{j=1}^m \varPi( {i\,a_1^\eta\,j} )(1-p(j \big\vert {i\,a_1^\eta} )) \leq 1
\end{equation*}
\end{linenomath*}

\begin{linenomath*}
\begin{equation*}
P \left ( \limsup_{n \to \infty} 
\frac{ \hat{\varDelta}_2(P_{\textbf{O}_n^\alpha} \Vert P_{\textbf{E}_n^\alpha})}{2\,\log(log(n))} \leq {\cal L} \right )=1.
\end{equation*}
\end{linenomath*}

\vspace{15pt}
\textit{\textbf{Case II:}} \,\, $\eta = \kappa - 1. $

In accordance with the following

\begin{linenomath*}
\begin{eqnarray*}
\lim_{n \to \infty} \frac{N(a_1^\eta\,\vert\,\textbf{X}_1^n)}{n} &=& 
\lim_{n \to \infty} \sum_{a \in E} \frac{N(a\,a_1^\eta\,\vert\,\textbf{X}_1^n)}{n} =
\sum_{a \in E} \varPi(a\,a_1^\eta) \,\,\,\, a.s.\\
\lim_{n \to \infty} \frac{N(i\,a_1^\eta\,\vert\,\textbf{X}_1^n)}{n} &=& \varPi(i\,a_1^\eta) \,\,\,\, a.s.\\
\end{eqnarray*}
\end{linenomath*}

we can obtain, as in previous case

\begin{linenomath*}
\begin{eqnarray*}
\lim_{n \to \infty}\textbf{P}_{\textbf{E}_n^\alpha}(i,j) &=& 
\lim_{n \to \infty} \frac{{\textbf{E}_n^\alpha}(i,j)}{N(a_1^\eta\,\vert\,\textbf{X}_1^n)}=
\lim_{n \to \infty}\Big (\frac{N(i\,a_1^\eta\,\vert\,\textbf{X}_1^n)}{N(a_1^\eta\,\vert\,\textbf{X}_1^n)}
                     \frac{N(a_1^\eta\,j\,\vert\,\textbf{X}_1^n)}{N(a_1^\eta\,\vert\,\textbf{X}_1^n)} \Big )=\\
&=&  \frac{\varPi(i\,a_1^\eta)}{\sum_{a \in E} \varPi(a\,a_1^\eta)}\frac{\varPi(a_1^\eta\,j)}{\sum_{a \in E} \varPi(a\,a_1^\eta)} \neq 0,\\
\label{prob_expected_n_alpha_1}
\end{eqnarray*}
\end{linenomath*}

and

\begin{linenomath*}
\begin{eqnarray*}
\lim_{n \to \infty}\textbf{P}_{\textbf{O}_n^\alpha}(i,j) 
&=& \lim_{n \to \infty} \frac{{\textbf{O}_n^\alpha}(i,j)}{N(a_1^\eta\,\vert\,\textbf{X}_1^n)}=
\lim_{n \to \infty} \Big( \frac{N(i\,a_1^\eta\,j\,\vert\,\textbf{X}_1^n)}{N(a_1^\eta\,\vert\,\textbf{X}_1^n)} \Big)=\\
&=& \frac{\varPi(i\,a_1^\eta)\,p(j\,\vert\,i\,a_1^\eta)}{\sum_{a \in E} \varPi(a\,a_1^\eta)} \neq 0.
\label{prob_observed_n_alpha_2}
\end{eqnarray*}
\end{linenomath*}

\vspace{15pt}
Clearly, if  $\eta = \kappa - 1$, there exist  $\alpha = a_1^\eta \,\, {\scriptstyle \&} \,\,\, i,j \in E$ so that 

\begin{linenomath*}
\begin{equation*}
\lim_{n \to \infty}(\textbf{P}_{\textbf{O}_n^\alpha}(i,j) - \textbf{P}_{\textbf{E}_n^\alpha}(i,j)) \neq 0 
\end{equation*}
\end{linenomath*}
 
since, otherwise, it should imply that

\begin{linenomath*}
\begin{equation*}
 p(j\,\vert\,i\,a_1^\eta) = \frac{\varPi(a_1^\eta\,j)}{\sum_{a \in E} \varPi(a\,a_1^\eta)}
\end{equation*}
\end{linenomath*}

i.e. $p(j\,\vert\,i\,a_1^\eta)$ does not depend on $i \in E$, contradicting the assumption that the order $\kappa > \eta$.

\begin{linenomath*}
\begin{equation*}
 P \left ( \hat{\varDelta}_2 (P_{\textbf{O}_n^\alpha} \Vert P_{\textbf{E}_n^\alpha} ) = n\,O(1) \right ) = 1
\end{equation*}
\end{linenomath*}

and (\ref{hypo_false_conv}) is proved. $ \quad \checkmark$

\vspace{30pt}
\subsection{Local and Global Dependency Level}

Herein we define the Local Dependency Level and the Global Dependency Level.
\vspace{10pt}
\begin{definition} 
Let $\textbf{X}_n=\{X_i\}_{i=1}^n$ be a sample of a Markov chain $\mathbb{X}$ of order $\kappa \geq 0$
and $\hat{\varDelta}_2({\textit{P}}_{{\textbf{O}}_n^\alpha} \Vert {\textit{P}}_{{\textbf{E}}_n^\alpha})$ with $\alpha=a_1^\eta,
\,\,\, \eta \geq 0 \,\,$ as previously defined.
 
Let us assume that $V$ is  a $\chi^2$ random variable with $(m-1)^2$ degrees of freedom
where  ${\cal{P}}$ is the continuous strictly decreasing function ${\cal{P}}:{\mathbb{R}}^{+} \longrightarrow [0,1]$ 

\begin{linenomath*}
\begin{equation*}
{\cal{P}}(x) = P( V \geq x ), \,\,\, x \in {\mathbb{R}}^{+}.
\end{equation*}
\end{linenomath*}

We define the Local Dependency Level $\widehat{LDL}_n(a_1^\eta)$, for $\alpha=a_1^\eta$ as

\begin{linenomath*}
\begin{equation*}
\boxed{
 \widehat{LDL}_n(a_1^\eta) =   
\frac{\hat{\varDelta}_2({\textit{P}}_{{\textbf{O}}_n^\alpha} \Vert {\textit{P}}_{{\textbf{E}}_n^\alpha})}{2 \,\log(\log(n))},
} 
\end{equation*}
\end{linenomath*}

\vspace{10pt}
and the Global Dependency Level $\widehat{GDL}_n(\eta)$ as 

\begin{linenomath*}
\begin{equation*}
\boxed{
\widehat{GDL}_n(\eta) = 
{\cal{P}} \left (\sum_{a_1^\eta \, \in E^\eta}\,\Big (\frac{N(a_1^\eta\,\vert\,\textbf{X}_1^n)}{n}\Big )\,\widehat{LDL}_n(a_1^\eta) \right ).}
\quad {\scriptstyle \blacklozenge}
\end{equation*}
\end{linenomath*}
\end{definition}

\vspace{30pt}
Observe that, if the hypothesis $\textbf{H}_0^\alpha$ is true, then \,
$\forall a_1^\eta, \,\, \eta~\geq~\kappa,$

\begin{linenomath*}
\begin{equation}
P \left (  \liminf_{n \to \infty} \left ( \widehat{GDL}_n(\eta) \right ) \geq {\cal{P}}({\cal{L}})  \right ) = 1 
\label{GDL_conv_1}
\end{equation}
\end{linenomath*}

and for $\eta = \kappa -1$

\begin{linenomath*}
\begin{equation}
P \left ( \lim_{n \to \infty} \left (  \widehat{GDL}_n(\eta)  \right ) = {\cal{P}}(\infty)= 0  \right ) = 1.
\label{GDL_conv_2}
\end{equation}
\end{linenomath*}

\vspace{15pt}
By (\ref{GDL_conv_1}) and (\ref{GDL_conv_2}) it is clear that, for $n$ sufficiently large,

\begin{linenomath*}
\begin{equation*}
P \left (  \widehat{GDL}_n(\eta) \approx 0 \right ) = 1, \quad \eta =  \kappa - 1, 
\label{GDL_conv_3}
\end{equation*}
\end{linenomath*}

and

\begin{linenomath*}
\begin{equation*}
P \left (  \widehat{GDL}_n(\eta) \approx {\cal{P}}({\cal{L}}) \right ) = 1, \quad  \eta \geq \kappa. 
\end{equation*}
\end{linenomath*}

and consequently, for a multiple stationary Markov chain $\mathbb{X}_{n \geq 1}$ of order $\kappa$

\begin{linenomath*}
\begin{eqnarray*}
\kappa &=& 0 \,\,\, \Leftrightarrow \,\,\, \lim_{n \to \infty} \widehat{GDL}_n(\eta) = {\cal P(L)}, \,\, \eta=0,1,..,B\,,\\
\kappa &=&  \max_{0\,\leq\,\eta\,\leq\,B} \Big\{\,\, \eta \, : \, \lim_{n \to \infty} \widehat{GDL}_n(\eta) = 0  \,\, \Big \} + 1.
\end{eqnarray*}
\end{linenomath*}

\vspace{40pt}
Finally, let us define the Markov chain order estimator based on the information contained in the vector $GDL_n.$

\begin{definition}
Given a fixed number $0 < B \in \mathbb{N}$, let us define the set ${\cal S} = \{0,1\}^{B+1}$ and the application 
$\,\, T\,:\,{\cal S} \rightarrow \mathbb{N}$

\begin{linenomath*}
\begin{eqnarray*}
T(s) &=& -1 \,\,\, \Leftrightarrow \,\,\, s_i=1, \,\,\, i=0,1,..,B \\
T(s) &=& \max_{0 \leq i \leq B} \left \{i \, : \, s_i=0, \, s_{i+1} = {\cal P(L)} \right \}, \,\,\, s=(s_0,s_1,...,s_B).
\quad {\scriptstyle \blacklozenge}
\end{eqnarray*}
\end{linenomath*}
\end{definition}

\vspace{25pt}

\begin{definition} 
Let $\textbf{X}_n=\{X_i\}_{i=1}^n$ be 
a sample  for the Markov chain $\mathbb{X}$ of order $\kappa$, $0 \leq \kappa < B \in \mathbb{N}$  and  $\{\widehat{GDL}_n(i) \}_{i=1}^B$ as above. 
We define the order's estimator $\widehat{\kappa}_{GDL}(\textbf{X}_n)$ as 

\begin{linenomath*}
\begin{equation*}
\widehat{\kappa}_{GDL}(\textbf{X}_n) = T( \sigma_n ) + 1
\end{equation*}
\end{linenomath*}

with $\sigma_n \in {\cal S}$ so that $\,\, \forall \, s \in {\cal S}$

\begin{linenomath*}
\begin{equation*}
\sum_{i=0}^B \big (\widehat{GDL}_n(i) - \sigma_n(i) \big)^2 \leq  \sum_{i=0}^B (\widehat{GDL}_n(i) - s(i))^2. 
\quad {\scriptstyle \blacklozenge}
\end{equation*}
\end{linenomath*}
\end{definition}

\vspace{20pt}
By (\ref{GDL_conv_1}),(\ref{GDL_conv_2}) and (\ref{GDL_conv_3}) it is clear that, for $n$ large enough,  
$\{GDL_n(i) \}_{i=1}^B$ satisfies the hypothesis of  
therefore, the order estimator converges almost surely to its value, i.e.,

\begin{linenomath*}
\begin{equation}
\boxed{\quad
P \left ( \lim_{n \to \infty}  \widehat{\kappa}_{GDL}(\textbf{X}_n) = \kappa \right )=1, \,\,\,\ \kappa = 0,1,2,..,B.\quad}
\label{order_conv_1} 
\end{equation}
\end{linenomath*}

\vspace{30pt}
\section{Numerical Simulations}
In what follows we shall compare the non-asymptotic performance, mainly for small samples, of some of the most used
Markov chains order estimators.
Recalling the previous notations 
$\alpha = (a_1,...,a_{k+1})= a_1^{k+1}$, \, ${N(i\,a_1^{k+1}\,\vert\,\textbf{X}_1^n)}$ \, as in (\ref{number_of_strings})
and denoting

\begin{linenomath*}
\begin{equation*}
{\hat L}(\eta) = \varPi_{a_1^{\eta+1}}
{\left [ \frac{{N(i\,a_1^{\eta+1}\,\vert\,\textbf{X}_1^n)}}{{N(i\,a_1^\eta \,\vert\,\textbf{X}_1^n)}} \right ]}
                                                                                          ^{N(i\,a_1^{\eta+1}\,\vert\,\textbf{X}_1^n)}
\end{equation*}
\end{linenomath*}

the estimators of the Markov chain order $\kappa$, are defined, under the hypothesis:
\begin{linenomath*}
\begin{equation*}
\boxed{ 
\textit{There exist a known}\,\, B \,\, \textit{so that} \,\, 0 \, \leq \kappa \,\, \leq B}
\end{equation*}
\end{linenomath*}
 
 as
\begin{linenomath*}
\begin{eqnarray*}
 &&{\widehat \kappa}_{AIC} = \text{argmin}\{ AIC(\eta)\,;\, \eta=0,1,...,B\}, \\
 &&{\widehat \kappa}_{BIC} = \text{argmin}\{ BIC(\eta)\,;\, \eta=0,1,...,B\}, \\
 &&{\widehat \kappa}_{EDC} = \text{argmin}\{ EDC(\eta)\,;\, \eta=0,1,...,B\},
\end{eqnarray*}
\end{linenomath*}

where

\begin{linenomath*}
\begin{eqnarray*}
  AIC(\eta) = -2 \log {\hat L}(\eta) &+&  {\vert E \vert}^{\eta+1} \,\,2\,(\vert E \vert - 1), \\  
  BIC(\eta) = -2 \log {\hat L}(\eta) &+&  {\vert E \vert}^{\eta+1} \,\,2\,(\vert E \vert - 1) \left( \frac{\log(n)}{2} \right ), \\
  EDC(\eta) = -2 \log {\hat L}(\eta) &+&  {\vert E \vert}^{\eta+1} \,\,2\,(\vert E \vert - 1) \left(\frac{\log\log(n)}{2 (\vert E \vert - 1)}\right),\\
  AIC(\eta)  \,\, \leq   &EDC(\eta)&  \leq \,\, BIC(\eta).
\end{eqnarray*}
\end{linenomath*}

\vspace{20pt}
Clearly, for a given $\eta$, the order estimator $GDL(\eta)$, as well as $AIC(\eta)$ \cite{Tong},
$BIC(\eta)$ \cite{Schwarz} and $EDC(\eta)$ \cite{Zhao,Dorea} contain
much of the information concerning the sample's relative dependency, nevertheless  numerical simulations as well as 
theoretical considerations anticipates a great deal of variability for small samples.

The following numerical simulation, based on an algorithm due to Raftery\cite{Raftery}, starts on with the generation 
of a  Markov chain transition matrix,  $\textbf{Q}=(q_{i_1i_2...i_{\kappa};i_{\kappa +1}})$ with entries 

\begin{linenomath*}
\begin{equation}
q_{i_1i_2...i_{\kappa};i_{\kappa +1}} = \sum_{t=1}^\kappa \lambda_{i_t} R(i_{\kappa +1},i_t), \,\, 
1 \leq i_t,i_{\kappa + 1} \leq m.
\label{markov_matrix_generation}
\end{equation}
\end{linenomath*}

where the matrix

\begin{linenomath*}
\begin{equation*}
{R}(i,j), \,\, 0 \leq i,j \leq m, \quad \sum_{i=1}^m R(i,j)=1, \,\, 1 \leq j \leq m
\end{equation*}
\end{linenomath*}

and the positive numbers 

\begin{linenomath*}
\begin{equation*}
\{ \lambda_i \}_{i=1}^\kappa, \,\, \sum_{i=1}^\kappa \lambda_i =1
\end{equation*}
\end{linenomath*}

are arbitrarily chosen in advance.

Once the  matrix $\textbf{Q}=(q_{i_1i_2...i_{\kappa};i_{\kappa +1}})$ is obtained, two hundreds
replications of the  Markov chain sample of size $n$, space state $E$ and transition matrix $\textbf{Q}$ 
are generated to  compare $GDL(\eta)$ performance against the standards, well known and 
already established order estimators just mentioned above.

\vspace{10pt}
Katz(1981) \cite{Katz} obtained the asymptotic distribution of ${\widehat \kappa}_{AIC}$ and proved 
its inconsistency showing 
the existence of a positive probability to overestimate the order. See also Shibata(1976) \cite{Shibata}. 

On the contrary Schwarz (1978) \cite{Schwarz} and Zhao(2001) \cite{Zhao} proved strong consistency for the estimators 
${\widehat \kappa}_{BIC}$ and ${\widehat \kappa}_{EDC}$, respectively. 

It is quite intuitive that  the random information regarding the order of a Markov chain, is spread over
an exponentially growing set of empirical distributions $\Theta$ with $\vert \Theta \vert = m^{B+1}$,
where $\textbf{B}$ is the maximum integer $\eta$, as in $\alpha = (i_1i_2...i_\eta)$. It seems  reasonable to think that a 
small \textit{viable} sample, i.e. samples able to retrieve enough information to estimate the chain order, should have  size 
$\,\,n \approx O(m^{B+1}).$ Keeping in mind that for the present numerical simulation, the maximum length to 
be used is $B=5$, from now on the sample
sizes for $\vert E \vert = 3$ and $\vert E \vert = 4$ should be $n \approx 1.500$ and $n \approx 5.000$, respectively.

Finally, after applying all estimators to
each one of the replicated samples, the final results are registered in the form of tables.

\vspace{60pt}
{\textbf{Case I: \,\, Markov \,\, Chain \,\, Examples \,\, with }\,\,  $\kappa = 0, \,\, \vert E \vert = 3.$ }

\vspace{10pt}
Firstly, we choose the matrix $ \{Q_1, Q_2, Q_3\}$ to produce  samples with sizes $500 \leq n \leq 2.000$,
originated from Markov chains of order $\kappa=0$  with quite different probability distributions.

\begin{linenomath*}
\[ 
Q_{1}= 
\left[ 
\begin{array}{cccc}
\scriptstyle{0.33} & \scriptstyle{0.335} & \scriptstyle{0.335} \\
\scriptstyle{0.33} & \scriptstyle{0.335} & \scriptstyle{0.335} \\
\scriptstyle{0.33} & \scriptstyle{0.335} & \scriptstyle{0.335} 
\end{array} 
\right]
,\,
Q_{2}= 
\left[ 
\begin{array}{cccc}
\scriptstyle{0.05} & \scriptstyle{0.475} & \scriptstyle{0.475} \\
\scriptstyle{0.05} & \scriptstyle{0.475} & \scriptstyle{0.475} \\
\scriptstyle{0.05} & \scriptstyle{0.475} & \scriptstyle{0.475} 
\end{array} 
\right]
,\,
Q_{3}= 
\left[ 
\begin{array}{cccc}
\scriptstyle{0.05} & \scriptstyle{0.05} & \scriptstyle{0.90} \\
\scriptstyle{0.05} & \scriptstyle{0.05} & \scriptstyle{0.90} \\
\scriptstyle{0.05} & \scriptstyle{0.05} & \scriptstyle{0.90}  
\end{array} 
\right].
\]
\end{linenomath*}

\begin{tabular}{|c||c|c|c|c||c|c|c|c||c|c|c|c||c|c|c|c|}
                                                     \hline 
\multicolumn{13}{|c|}{\small 
     \textbf{\textit{ $\vert E \vert = 3 \qquad \leftrightarrow  \qquad \small  \kappa = 0  
                                         \qquad \leftrightarrow  \qquad \scriptstyle {\lambda_i\,=\,1/3, \,\, i\,=\,1,2,3.}$}}}\\
                                                     \hline
                                      &\multicolumn{4}{|c||}{\small $Q_1$}  
                                      &\multicolumn{4}{|c||}{\small $Q_1$}  
                                      &\multicolumn{4}{|c||}{\small $Q_1$}\\
                                                    \hline
                                      &\multicolumn{4}{|c||}{$\scriptstyle  n\,=\,500 $}
                                      &\multicolumn{4}{|c||}{$\scriptstyle  n\,=\,1.000 $}
                                      &\multicolumn{4}{|c||}{$\scriptstyle  n\,=\,1.500$ }\\
                                                    \hline \hline
$\scriptstyle k$     & $\scriptstyle Aic$ &  $\scriptstyle Bic$ & $\scriptstyle Edc$  & $\scriptstyle Gdl$ &
                       $\scriptstyle Aic$ &  $\scriptstyle Bic$ & $\scriptstyle Edc$  & $\scriptstyle Gdl$ &
                       $\scriptstyle Aic$ &  $\scriptstyle Bic$ & $\scriptstyle Edc$  & $\scriptstyle Gdl$ \\
                       	                             \hline \hline
$\scriptstyle 0$     & $\scriptstyle 75.5\%$ & $\scriptstyle 100\%$ & $\scriptstyle 100\% $ & $\scriptstyle 99\%   $ & 
                       $\scriptstyle 80\%  $ & $\scriptstyle 100\%$ & $\scriptstyle 100\% $ & $\scriptstyle 99.5\% $ & 
                       $\scriptstyle 71.5\%$ & $\scriptstyle 100\%$ & $\scriptstyle 100\% $ & $\scriptstyle 99\%   $ \\
                                                            \hline 
$\scriptstyle 1$     & $\scriptstyle 24.5\%$ & $\scriptstyle     $ & $\scriptstyle     $ & $\scriptstyle 1\%    $ & 
                       $\scriptstyle 18\%  $ & $\scriptstyle     $ & $\scriptstyle     $ & $\scriptstyle 0.5\%  $ & 
                       $\scriptstyle 22.5\%$ & $\scriptstyle     $ & $\scriptstyle     $ & $\scriptstyle 1\%    $ \\
                                                            \hline  
$\scriptstyle 2$     & $\scriptstyle      $ & $\scriptstyle  $ & $\scriptstyle    $ & $\scriptstyle    $ & 
                       $\scriptstyle  2\% $ & $\scriptstyle  $ & $\scriptstyle    $ & $\scriptstyle    $ & 
                       $\scriptstyle  6\% $ & $\scriptstyle  $ & $\scriptstyle    $ & $\scriptstyle    $ \\
                                                            \hline
$\scriptstyle 3$     & $\scriptstyle   $ & $\scriptstyle      $ & $\scriptstyle   $ & $\scriptstyle   $ & 
                       $\scriptstyle   $ & $\scriptstyle      $ & $\scriptstyle   $ & $\scriptstyle   $ & 
                       $\scriptstyle   $ & $\scriptstyle      $ & $\scriptstyle   $ & $\scriptstyle   $ \\
                                                            \hline
$\scriptstyle 4$     & $\scriptstyle \phantom{100\%}    $ & $\scriptstyle    $ & $\scriptstyle    $ & $\scriptstyle \phantom{100\%}    $ & 
                       $\scriptstyle \phantom{100\%}    $ & $\scriptstyle    $ & $\scriptstyle    $ & $\scriptstyle \phantom{100\%}    $ & 
                       $\scriptstyle \phantom{100\%}    $ & $\scriptstyle    $ & $\scriptstyle    $ & $\scriptstyle \phantom{100\%}    $ \\
                                                       \hline \hline
\end{tabular}

\vspace{10pt}
\begin{tabular}{|c||c|c|c|c||c|c|c|c||c|c|c|c||c|c|c|c|}
                                                     \hline 
\multicolumn{13}{|c|}{\small 
     \textbf{\textit{ $\vert E \vert = 3 \qquad \leftrightarrow  \qquad \small  \kappa = 0  
                                         \qquad \leftrightarrow  \qquad \scriptstyle {\lambda_i\,=\,1/3, \,\, i\,=\,1,2,3.}$}}}\\
                                                     \hline
                                      &\multicolumn{4}{|c||}{\small $Q_2$}  
                                      &\multicolumn{4}{|c||}{\small $Q_2$}  
                                      &\multicolumn{4}{|c||}{\small $Q_2$}\\
                                                     \hline
                                      &\multicolumn{4}{|c||}{$\scriptstyle  n\,=\,1.000 $}
                                      &\multicolumn{4}{|c||}{$\scriptstyle  n\,=\,1.500 $}
                                      &\multicolumn{4}{|c||}{$\scriptstyle  n\,=\,500   $}\\
                                                    \hline \hline
$\scriptstyle k$     & $\scriptstyle Aic$ &  $\scriptstyle Bic$ & $\scriptstyle Edc$  & $\scriptstyle Gdl$ &
                       $\scriptstyle Aic$ &  $\scriptstyle Bic$ & $\scriptstyle Edc$  & $\scriptstyle Gdl$ &
                       $\scriptstyle Aic$ &  $\scriptstyle Bic$ & $\scriptstyle Edc$  & $\scriptstyle Gdl$ \\
                       	                             \hline \hline
$\scriptstyle 0$     & $\scriptstyle 63.5\%$ & $\scriptstyle 100\% $ & $\scriptstyle 100\% $ & $\scriptstyle  99\% $ & 
                       $\scriptstyle 63\%  $ & $\scriptstyle 100\% $ & $\scriptstyle 100\% $ & $\scriptstyle  99\%   $ & 
                       $\scriptstyle 59\%  $ & $\scriptstyle 100\% $ & $\scriptstyle 100\% $ & $\scriptstyle  99\% $ \\
                                                            \hline
$\scriptstyle 1$     & $\scriptstyle 29\%  $ & $\scriptstyle     $ & $\scriptstyle     $ & $\scriptstyle  1\%  $ & 
                       $\scriptstyle 34.5\%$ & $\scriptstyle     $ & $\scriptstyle     $ & $\scriptstyle  1\%  $ & 
                       $\scriptstyle 37\%  $ & $\scriptstyle     $ & $\scriptstyle     $ & $\scriptstyle  1\%  $\\
                                                            \hline
$\scriptstyle 2$     & $\scriptstyle 7.5\% $ & $\scriptstyle        $ & $\scriptstyle       $ & $\scriptstyle    $ & 
                       $\scriptstyle 2.5\% $ & $\scriptstyle        $ & $\scriptstyle       $ & $\scriptstyle    $ & 
                       $\scriptstyle 4\%   $ & $\scriptstyle        $ & $\scriptstyle       $ & $\scriptstyle    $ \\
                                                            \hline
$\scriptstyle 3$     & $\scriptstyle      $ & $\scriptstyle       $ & $\scriptstyle       $ & $\scriptstyle      $ & 
                       $\scriptstyle      $ & $\scriptstyle       $ & $\scriptstyle       $ & $\scriptstyle      $ & 
                       $\scriptstyle      $ & $\scriptstyle       $ & $\scriptstyle       $ & $\scriptstyle      $ \\
                                                            \hline
$\scriptstyle 4$     & $\scriptstyle \phantom{100\%}    $ & $\scriptstyle    $ & $\scriptstyle    $ & $\scriptstyle \phantom{100\%}    $ & 
                       $\scriptstyle \phantom{100\%}    $ & $\scriptstyle    $ & $\scriptstyle    $ & $\scriptstyle \phantom{100\%}    $ & 
                       $\scriptstyle \phantom{100\%}    $ & $\scriptstyle    $ & $\scriptstyle    $ & $\scriptstyle \phantom{100\%}    $ \\
                                                       \hline \hline
\end{tabular}

\vspace{10pt}
\begin{tabular}{|c||c|c|c|c||c|c|c|c||c|c|c|c||c|c|c|c|}
                \hline 
\multicolumn{13}{|c|}{\small 
     \textbf{\textit{ $\vert E \vert = 3 \qquad \leftrightarrow  \qquad \small  \kappa = 0  
                                         \qquad \leftrightarrow  \qquad \scriptstyle {\lambda_i\,=\,1/3, \,\, i\,=\,1,2,3.}$}}}\\
                                                     \hline
                                      &\multicolumn{4}{|c||}{\small $Q_3$}  
                                      &\multicolumn{4}{|c||}{\small $Q_3$}  
                                      &\multicolumn{4}{|c||}{\small $Q_3$}\\
                                                     \hline
                                      &\multicolumn{4}{|c||}{$\scriptstyle  n\,=\,1.000 $}
                                      &\multicolumn{4}{|c||}{$\scriptstyle  n\,=\,1.500 $}
                                      &\multicolumn{4}{|c||}{$\scriptstyle  n\,=\,2.000$ }\\
                                                    \hline \hline
$\scriptstyle k$     & $\scriptstyle Aic$ &  $\scriptstyle Bic$ & $\scriptstyle Edc$  & $\scriptstyle Gdl$ &
                       $\scriptstyle Aic$ &  $\scriptstyle Bic$ & $\scriptstyle Edc$  & $\scriptstyle Gdl$ &
                       $\scriptstyle Aic$ &  $\scriptstyle Bic$ & $\scriptstyle Edc$  & $\scriptstyle Gdl$ \\
                       	                             \hline \hline
$\scriptstyle 0$     & $\scriptstyle 43\%$ & $\scriptstyle 100\%$ & $\scriptstyle 100\% $ & $\scriptstyle 98\%$ & 
                       $\scriptstyle 47\%$ & $\scriptstyle 100\%$ & $\scriptstyle 99.5\%$ & $\scriptstyle 96\%$ & 
                       $\scriptstyle 46\%$ & $\scriptstyle 100\%$ & $\scriptstyle 100\% $ & $\scriptstyle 97\%$ \\
                                                            \hline
$\scriptstyle 1$     & $\scriptstyle 53\%  $ & $\scriptstyle     $ & $\scriptstyle      $ & $\scriptstyle 2\%  $ & 
                       $\scriptstyle 51.5\%$ & $\scriptstyle     $ & $\scriptstyle 0.5\%$ & $\scriptstyle 4\%  $ & 
                       $\scriptstyle 50.5\%$ & $\scriptstyle     $ & $\scriptstyle      $ & $\scriptstyle 2\%  $ \\
                                                            \hline
$\scriptstyle 2$     & $\scriptstyle  4\% $ & $\scriptstyle     $ & $\scriptstyle    $ & $\scriptstyle      $ & 
                       $\scriptstyle 1.5\%$ & $\scriptstyle     $ & $\scriptstyle    $ & $\scriptstyle      $ & 
                       $\scriptstyle 3.5\%$ & $\scriptstyle     $ & $\scriptstyle    $ & $\scriptstyle  1\% $ \\
                                                            \hline
$\scriptstyle 3$     & $\scriptstyle      $ & $\scriptstyle       $ & $\scriptstyle     $ & $\scriptstyle        $ & 
                       $\scriptstyle      $ & $\scriptstyle       $ & $\scriptstyle     $ & $\scriptstyle        $ & 
                       $\scriptstyle      $ & $\scriptstyle       $ & $\scriptstyle     $ & $\scriptstyle        $ \\
                                                            \hline
$\scriptstyle 4$     & $\scriptstyle \phantom{100\%} $ & $\scriptstyle    $ & $\scriptstyle    $ & $\scriptstyle \phantom{100\%}    $ & 
                       $\scriptstyle \phantom{100\%} $ & $\scriptstyle    $ & $\scriptstyle    $ & $\scriptstyle \phantom{100\%}    $ & 
                       $\scriptstyle \phantom{100\%} $ & $\scriptstyle    $ & $\scriptstyle    $ & $\scriptstyle \phantom{100\%}    $ \\
                                                       \hline \hline
\end{tabular}

\vspace{25pt}
Notice that for a fixed sample size $n = \{ 500, 1.000, 1.500, 2.000 \}$, the order estimator $\widehat{\kappa}_{AIC}$ 
steadily overestimate the real order $\kappa = 0$ with the excessiveness depending on the  probability distribution of the Markov chain. 
Differently, the order estimators $\widehat{\kappa}_{BIC}$, $\widehat{\kappa}_{EDC}$ and $\widehat{\kappa}_{GDL}$ show consistent performance, 
mainly obtaining the right order, free from the influence of the sample size and the generating matrix.
Regarding $\widehat{\kappa}_{BIC}$ and  $\widehat{\kappa}_{EDC}$ improved effect, most likely depends on their correcting factor,
$\frac{\log(n)}{2}$ and $\left( \frac{\log\log(n)}{2 (\vert E \vert - 1)} \right )$
which tend to decrease the estimated order.

\vspace{60pt}
{\textbf{Case II: \,\, Markov \,\, Chain \,\, Examples \,\, with }\,\,  $\kappa = 3, \vert E \vert = 3$ and 
$\kappa = \{2,3,0\}, \vert E \vert = 4$}

\vspace{10pt}
Secondly, we choose the matrix $ \{Q_4, Q_5\}$ to produce samples  with sizes 
$n \in \{ 500, 1.000, 1.500, 2.000 \}$,
originated from Markov chains for $\vert E\vert=3$ of order $\kappa=3$.

\vspace{30pt}
\begin{linenomath*}
\[ 
Q_{4}= 
\left[ 
\begin{array}{cccc} 
\scriptstyle{0.05} & \scriptstyle{0.05} & \scriptstyle{0.90} \\
\scriptstyle{0.05} & \scriptstyle{0.90} & \scriptstyle{0.05} \\
\scriptstyle{0.90} & \scriptstyle{0.05} & \scriptstyle{0.05} 
\end{array} 
\right]
,\qquad
Q_{5}= 
\left[ 
\begin{array}{cccc}
\scriptstyle{0.475} & \scriptstyle{0.475} & \scriptstyle{0.05} \\
\scriptstyle{0.475} & \scriptstyle{0.05} & \scriptstyle{0.475} \\
\scriptstyle{0.05} & \scriptstyle{0.475} & \scriptstyle{0.475} 
\end{array} 
\right].
\]
\end{linenomath*}

\vspace{20pt}
\begin{tabular}{|c||c|c|c|c||c|c|c|c||c|c|c|c||c|c|c|c|}
                \hline 
\multicolumn{13}{|c|}{\small 
     \textbf{\textit{ $\vert E \vert = 3 \qquad \leftrightarrow  \qquad \small  \kappa = 3  
                     \qquad \leftrightarrow  \qquad \scriptstyle \scriptstyle {\lambda_i\,=\,1/3, \,\, i\,=\,1,2,3.}$}}}\\
                                                     \hline
                                      &\multicolumn{4}{|c||}{\small $Q_4$}  
                                      &\multicolumn{4}{|c||}{\small $Q_4$}  
                                      &\multicolumn{4}{|c||}{\small $Q_4$}\\
                                                     \hline
                                      &\multicolumn{4}{|c||}{$\scriptstyle  n\,=\,1.000 $}
                                      &\multicolumn{4}{|c||}{$\scriptstyle  n\,=\,1.500 $}
                                      &\multicolumn{4}{|c||}{$\scriptstyle  n\,=\,2.000 $}\\
                                                    \hline \hline
$\scriptstyle k$     & $\scriptstyle Aic$ &  $\scriptstyle Bic$ & $\scriptstyle Edc$  & $\scriptstyle Gdl$ &
                       $\scriptstyle Aic$ &  $\scriptstyle Bic$ & $\scriptstyle Edc$  & $\scriptstyle Gdl$ &
                       $\scriptstyle Aic$ &  $\scriptstyle Bic$ & $\scriptstyle Edc$  & $\scriptstyle Gdl$ \\
                       	                             \hline \hline
$\scriptstyle 0$     & $\scriptstyle \phantom{100\%}  $ & $\scriptstyle \phantom{100\%}  $ & $\scriptstyle  $ & $\scriptstyle \phantom{100\%}  $ & 
                       $\scriptstyle \phantom{100\%}  $ & $\scriptstyle \phantom{100\%}  $ & $\scriptstyle  $ & $\scriptstyle \phantom{100\%}  $ & 
                       $\scriptstyle \phantom{100\%}  $ & $\scriptstyle \phantom{100\%}  $ & $\scriptstyle  $ & $\scriptstyle  \phantom{100\%} $ \\
                                                            \hline 
$\scriptstyle 1$     & $\scriptstyle   $ & $\scriptstyle     $ & $\scriptstyle     $ & $\scriptstyle     $ & 
                       $\scriptstyle   $ & $\scriptstyle     $ & $\scriptstyle     $ & $\scriptstyle     $ & 
                       $\scriptstyle   $ & $\scriptstyle     $ & $\scriptstyle     $ & $\scriptstyle     $ \\
                                                            \hline  
$\scriptstyle 2$     & $\scriptstyle      $ & $\scriptstyle 99.5\% $ & $\scriptstyle  88.5\%  $ & $\scriptstyle  41\%  $ & 
                       $\scriptstyle      $ & $\scriptstyle 76.5\% $ & $\scriptstyle  16.5\%  $ & $\scriptstyle   5\%  $ & 
                       $\scriptstyle      $ & $\scriptstyle 17\%   $ & $\scriptstyle  0.5\%   $ & $\scriptstyle  1\%   $ \\
                                                            \hline
$\scriptstyle 3$     & $\scriptstyle 100\%  $ & $\scriptstyle 0.5\%   $ & $\scriptstyle 11.5\% $ & $\scriptstyle 59\% $ & 
                       $\scriptstyle 100\%  $ & $\scriptstyle 23.5\%  $ & $\scriptstyle 83.5\% $ & $\scriptstyle 95\% $ & 
                       $\scriptstyle 100\%  $ & $\scriptstyle  83\%   $ & $\scriptstyle 99.5\% $ & $\scriptstyle 99\% $ \\
                                                            \hline
$\scriptstyle 4$     & $\scriptstyle    $ & $\scriptstyle    $ & $\scriptstyle    $ & $\scriptstyle    $ & 
                       $\scriptstyle    $ & $\scriptstyle    $ & $\scriptstyle    $ & $\scriptstyle    $ & 
                       $\scriptstyle    $ & $\scriptstyle    $ & $\scriptstyle    $ & $\scriptstyle    $ \\
                                                       \hline \hline
\end{tabular}

\vspace{30pt}
\begin{tabular}{|c||c|c|c|c||c|c|c|c||c|c|c|c||c|c|c|c|}
                \hline 
\multicolumn{13}{|c|}{\small \textbf{\textit{ $\vert E \vert = 3 \quad \leftrightarrow  \quad \small  \kappa = 3
               \quad \leftrightarrow  \quad \scriptstyle \scriptstyle {\lambda_i\,=\,1/3, \,\, i\,=\,1,2,3.}$}}}\\
                                                     \hline
                                        &\multicolumn{4}{|c||}{\small $Q_5$}
                                        &\multicolumn{4}{|c||}{\small $Q_5$}
                                        &\multicolumn{4}{|c||}{\small $Q_5$}\\
                                                     \hline
                                      &\multicolumn{4}{|c||}{$\scriptstyle  n\,=\,1.000 $}
                                      &\multicolumn{4}{|c||}{$\scriptstyle  n\,=\,1.500 $}
                                      &\multicolumn{4}{|c||}{$\scriptstyle  n\,=\,2.500$ }\\
                                                    \hline \hline
$\scriptstyle k$     & $\scriptstyle Aic$ &  $\scriptstyle Bic$ & $\scriptstyle Edc$  & $\scriptstyle Gdl$ &
                       $\scriptstyle Aic$ &  $\scriptstyle Bic$ & $\scriptstyle Edc$  & $\scriptstyle Gdl$ &
                       $\scriptstyle Aic$ &  $\scriptstyle Bic$ & $\scriptstyle Edc$  & $\scriptstyle Gdl$ \\
                       	                             \hline \hline
$\scriptstyle 0$     & $\scriptstyle       $ & $\scriptstyle 0.5\% $ & $\scriptstyle         $ & $\scriptstyle         $ & 
                       $\scriptstyle       $ & $\scriptstyle       $ & $\scriptstyle         $ & $\scriptstyle         $ & 
                       $\scriptstyle       $ & $\scriptstyle       $ & $\scriptstyle         $ & $\scriptstyle         $ \\
                                                            \hline
$\scriptstyle 1$     & $\scriptstyle     $ & $\scriptstyle 92.5\% $ & $\scriptstyle 69.5\% $ & $\scriptstyle  6.5\%  $ & 
                       $\scriptstyle     $ & $\scriptstyle 54.5\% $ & $\scriptstyle 19.5\% $ & $\scriptstyle  1\%    $ & 
                       $\scriptstyle     $ & $\scriptstyle        $ & $\scriptstyle     $ & $\scriptstyle            $ \\
                                                            \hline
$\scriptstyle 2$     & $\scriptstyle 16.5\%$ & $\scriptstyle 7\%   $ & $\scriptstyle  30.5\% $ & $\scriptstyle 92\%  $ & 
                       $\scriptstyle  2\%  $ & $\scriptstyle 45.5\%$ & $\scriptstyle  80.5\% $ & $\scriptstyle 80.5\%$ & 
                       $\scriptstyle       $ & $\scriptstyle 100\% $ & $\scriptstyle  98.5\% $ & $\scriptstyle 8.5\% $ \\
                                                            \hline
$\scriptstyle 3$     & $\scriptstyle 83.5\%  $ & $\scriptstyle      $ & $\scriptstyle       $ & $\scriptstyle 1.5\% $ & 
                       $\scriptstyle 98\%    $ & $\scriptstyle      $ & $\scriptstyle       $ & $\scriptstyle 18.5\%$ & 
                       $\scriptstyle 100\%   $ & $\scriptstyle      $ & $\scriptstyle 1.5\% $ & $\scriptstyle 91.5\%$ \\
                                                            \hline
$\scriptstyle 4$     & $\scriptstyle    $ & $\scriptstyle    $ & $\scriptstyle    $ & $\scriptstyle    $ & 
                       $\scriptstyle    $ & $\scriptstyle    $ & $\scriptstyle    $ & $\scriptstyle    $ & 
                       $\scriptstyle    $ & $\scriptstyle    $ & $\scriptstyle    $ & $\scriptstyle    $ \\
                                                       \hline \hline
\end{tabular}

\vspace{40pt}
For $\vert E \vert = 3$\, , \, $\kappa = 3$ the estimator $\widehat{\kappa}_{AIC}$ overestimate the order 
in a lesser extent than the previous case, while $\widehat{\kappa}_{BIC}$ and $\widehat{\kappa}_{EDC}$ overweighted  
by the respective constants   $\frac{\log(n)}{2}$ and $\left( \frac{\log\log(n)}{2 (\vert E \vert - 1)} \right )$, 
underestimate the order more than it was supposed to be. 
Concerning $\widehat{\kappa}_{GDL}$, it rapidly converges to the right order depending on the sample size $n$.

For $\vert E \vert = 4$ the greater complexity of a Markov chain of order $\kappa = 3$ impose the use of larger sample size
for estimators to acomplish some reliability.
Finally, we choose the matrix $ \{ Q_6, Q_7 \}$ to produce samples with size $n = 5.000$,
originated from Markov chains of order $\kappa \in \{ 2,3,0 \}$  like in the previous cases.

\vspace{30pt}
\begin{linenomath*}
\[ 
Q_{6}= 
\left[ 
\begin{array}{cccc}
\scriptstyle{0.05} & \scriptstyle{0.05} & \scriptstyle{0.05} & \scriptstyle{0.85} \\
\scriptstyle{0.05} & \scriptstyle{0.05} & \scriptstyle{0.85} & \scriptstyle{0.05} \\
\scriptstyle{0.05} & \scriptstyle{0.85} & \scriptstyle{0.05} & \scriptstyle{0.05} \\
\scriptstyle{0.85} & \scriptstyle{0.05} & \scriptstyle{0.05} & \scriptstyle{0.05}
\end{array} 
\right]
,\qquad
Q_{7}= 
\left[ 
\begin{array}{cccc}
\scriptstyle{0.05} & \scriptstyle{0.05} & \scriptstyle{0.05} & \scriptstyle{0.85} \\
\scriptstyle{0.05} & \scriptstyle{0.05} & \scriptstyle{0.05} & \scriptstyle{0.85} \\
\scriptstyle{0.05} & \scriptstyle{0.05} & \scriptstyle{0.05} & \scriptstyle{0.85} \\
\scriptstyle{0.05} & \scriptstyle{0.05} & \scriptstyle{0.05} & \scriptstyle{0.85}
\end{array} 
\right].
\]
\end{linenomath*}

\begin{tabular}{|c||c|c|c|c||c|c|c|c||c|c|c|c||c|c|c|c|}
                \hline 
\multicolumn{13}{|c|}{\small \textbf{\textit{ $\vert E \vert \,=\, 4 \,\,\,\,\,\,\, \leftrightarrow \,\,\,\,\,\,\,\, n\,=\, 5.000   $}} } \\
                                                     \hline
     &\multicolumn{4}{|c||}{\small $Q_6  \Leftrightarrow \,\, \scriptstyle{ \lambda_i\,=\,1/2, \,\, i\,=\,1,2.  }$}
     &\multicolumn{4}{|c||}{\small $Q_6  \Leftrightarrow \,\, \scriptstyle{ \lambda_i\,=\,1/3, \,\, i\,=\,1,2,3.}$}
     &\multicolumn{4}{|c|}{\small  $Q_7  \Leftrightarrow \,\, \scriptstyle{ \lambda_i\,=\,1/3, \,\, i\,=\,1,2,3.}$}\\
                                                      \hline
                                      &\multicolumn{4}{|c||}{\small $ \kappa\,=\, 2 $}
                                      &\multicolumn{4}{|c||}{\small $ \kappa\,=\, 3 $}
                                      &\multicolumn{4}{|c|}{\small  $ \kappa\,=\, 0 $}\\
                                                    \hline \hline
$\scriptstyle k$     & $\scriptstyle Aic$ &  $\scriptstyle Bic$ & $\scriptstyle Edc$  & $\scriptstyle Gdl$ &
                       $\scriptstyle Aic$ &  $\scriptstyle Bic$ & $\scriptstyle Edc$  & $\scriptstyle Gdl$ &
                       $\scriptstyle Aic$ &  $\scriptstyle Bic$ & $\scriptstyle Edc$  & $\scriptstyle Gdl$ \\
                       	                             \hline \hline
$\scriptstyle 0$     & $\scriptstyle    $ & $\scriptstyle    $ & $\scriptstyle         $ & $\scriptstyle      $ & 
                       $\scriptstyle    $ & $\scriptstyle    $ & $\scriptstyle         $ & $\scriptstyle      $ & 
                       $\scriptstyle 85\%$ & $\scriptstyle 100\%$ & $\scriptstyle 100\%$ & $\scriptstyle 100\%$ \\
                                                            \hline
$\scriptstyle 1$     & $\scriptstyle     $ & $\scriptstyle    $ & $\scriptstyle    $ & $\scriptstyle    $ & 
                       $\scriptstyle     $ & $\scriptstyle    $ & $\scriptstyle    $ & $\scriptstyle    $ & 
                       $\scriptstyle 15\%$ & $\scriptstyle    $ & $\scriptstyle    $ & $\scriptstyle    $ \\
                                                            \hline
$\scriptstyle 2$     & $\scriptstyle 100\%$ & $\scriptstyle 100\% $ & $\scriptstyle 100\%$ & $\scriptstyle 100\%$ & 
                       $\scriptstyle      $ & $\scriptstyle   99\%$ & $\scriptstyle      $ & $\scriptstyle   4\%$ & 
                       $\scriptstyle      $ & $\scriptstyle       $ & $\scriptstyle      $ & $\scriptstyle      $ \\
                                                            \hline
$\scriptstyle 3$     & $\scriptstyle      $ & $\scriptstyle      $ & $\scriptstyle      $ & $\scriptstyle       $ & 
                       $\scriptstyle 100\%$ & $\scriptstyle   1\%$ & $\scriptstyle 100\%$ & $\scriptstyle   96\%$ & 
                       $\scriptstyle      $ & $\scriptstyle      $ & $\scriptstyle      $ & $\scriptstyle       $ \\
                                                            \hline
$\scriptstyle 4$     & $\scriptstyle    $ & $\scriptstyle    $ & $\scriptstyle    $ & $\scriptstyle    $ & 
                       $\scriptstyle    $ & $\scriptstyle    $ & $\scriptstyle    $ & $\scriptstyle    $ & 
                       $\scriptstyle    $ & $\scriptstyle    $ & $\scriptstyle    $ & $\scriptstyle    $ \\
                                                            \hline
$\scriptstyle 5$     & $\scriptstyle    $ & $\scriptstyle    $ & $\scriptstyle    $ & $\scriptstyle    $ & 
                       $\scriptstyle    $ & $\scriptstyle    $ & $\scriptstyle    $ & $\scriptstyle    $ & 
                       $\scriptstyle    $ & $\scriptstyle    $ & $\scriptstyle    $ & $\scriptstyle    $ \\
                                                            \hline
$\scriptstyle 6$     & $\scriptstyle    $ & $\scriptstyle    $ & $\scriptstyle    $ & $\scriptstyle    $ & 
                       $\scriptstyle    $ & $\scriptstyle    $ & $\scriptstyle    $ & $\scriptstyle    $ & 
                       $\scriptstyle    $ & $\scriptstyle    $ & $\scriptstyle    $ & $\scriptstyle    $ \\ 
                                                            \hline\hline
\end{tabular}

\vspace{25pt}
For the order for $\vert E \vert = 4$, $\kappa =0$, apparently  $\widehat{\kappa}_{AIC}$ keeps overestimating 
the order in some degree, while $\widehat{\kappa}_{BIC}$ as in example $\kappa = 3$
severely underestimate the order, presumably due to the excessive weight of the correcting factors $\frac{\log(n)}{2}$. 
On the contrary $\widehat{\kappa}_{EDC}$ and $\widehat{\kappa}_{GDL}$ behaves quite well in same setting.

\vspace{30pt}
\section{Conclusion}
The pioneer research started with the contributions of Bartlett\cite{Bartlett}, Hoel\cite{Hoel}, Good \cite{Good}, 
Anderson \& Goodman \cite{Anderson-Goodman}, Billingsley(\cite{Billingsley1}, \cite{Billingsley2}) among others, where they 
developed tests of hypothesis for the estimation of the order of a given Markov chain.

Later on these procedures were adapted and improved with the used of \textit{Penalty Functions} (Tong\cite{Tong}, Katz\cite{Katz}) 
together with other tools created in the realm of Models Selection (Akaike\cite{Akaike}, Schwarz\cite{Schwarz}).
Since then, there have been a considerable  number of subsequent contributions on this subject, several of them consisting in the
enhancement of the already  existing  techniques (Csiszar\cite{CsiszarShields_1}, Zhao et all\cite{Zhao}).

In this notes we propose a new Markov chain order estimator based on a  different idea which makes it behave 
in a quite different form.
This estimator is strongly consistent and  more efficient than AIC (inconsistent),  outperforming 
the well established and consistent BIC and EDC, mainly on relatively small samples.

\end{document}